\documentclass{article}
\usepackage{hyperref, url}
\usepackage{amscd,amsmath,amsthm,amssymb}
\usepackage{enumerate,varioref}
\usepackage{epsfig}
\usepackage{graphicx}
\usepackage{subfig}
\usepackage{mathtools}
\usepackage{tikz}
\usepackage{algorithmicx}
\usepackage{algorithm}
\usepackage{algpseudocode}
\usepackage{xcolor}
\newtheorem{thm}{Theorem}

\newtheorem{lem}[thm]{Lemma}

\theoremstyle{definition}
\newtheorem{defn}[thm]{Definition}
\theoremstyle{remark}
\newtheorem{ex}[thm]{Example}

\newcommand{\Z}{\mathbb{Z}}

\title{A Note on Consistent Rotation Maps of Graph Cartesian Products}
\author{Clark Alexander\footnote{Chicago Quantum, Independent} \\ email: \href{mailto:gcalexander1981@gmail.com}{the author}}

\begin{document}
	
	\maketitle
	
	\begin{abstract}
	Given two regular graphs with consistent rotation maps, we produce a constructive method for a consistent rotation map on their Cartesian product.  This method will be given as a simple set of rules of addition and table look ups.  We assume that the combinatorial construction of both consistent rotation maps has occurred before we construct the Cartesian product.
	\end{abstract}

	\tableofcontents
	
	\section{Introduction}
	
	The motivation for this work comes from constructing an efficient method for discrete time quantum walks (DTQW) on regular graphs.  Given that the variance of a quantum random walk can diverge wildly from that of a classical random walk \cite{A, AKR, K} the original aim of this work was to construct DTQW
	on regular expander graphs.  This led us to \cite{zig-zag} and the definition of the rotation map.  Further, in order to implement a coined DTQW rather than the style of Szegedy, Kendon, Ambainis, et al. \cite{K} the shifting operator requires a consistent rotation map to insure the overall operator is unitary.  In a greedy algorithm which builds a rotation map, one gets a partially hyperbolic system.
	
	Indeed in future works we will construct coined DTQW using consistent rotation maps.  At present we wish to build up the body of literature on rotation maps, as at present the literature is woefully thin.
	
	The rest of this work flows as follows: In \S2 we present some algebraic basics of rotation maps and Cartesian products of graphs. We then give the main method for building a consistent rotation map for a Cartesian product of graphs.  
	
	In \S3 we produce consistent rotation maps on several families of graphs including, Hypercubes, cycle graphs, complete graphs, complete bipartite graphs $(K_{n,n})$, and generalized Petersen graphs.  In some cases these rotation maps are so simple that they can be generated in two lines of code and other are so simple that we can simply write them down from a well-drawn picture.  We also show that these basic families of graphs with Cartesian products open up several other families for free e.g. square lattice tori.
	
	In \S4 we explore future work wherein we will apply these rotation maps to coined DTQW on regular graphs. We will also discuss some methods by which we can construct consistent rotation maps on arbitrary regular graphs.
	
	\section{The Method}
	\subsection{Defining Rotation Maps}
	
	In this section we'll transform the definition of rotation maps given in \cite{zig-zag} to something which we can compute in matrices.  We will show that these definitions are equivalent and then a method to read out the original rotation map formulation from our matrix formulation.
	
	\begin{defn}
		Given a regular graph $G = (V,E)$ with degree of regularity $d$, the rotation map on $G$ is a function
		\begin{equation}
		Rot_G : [|V|] \times [d] \rightarrow [|V|]\times [d]
		\end{equation}
		where $Rot_G(v,i) = (w,j)$ with the $i^{th}$ edge leaving from vertex $v$ enters vertex $w$, and the $j^{th}$ edge leaving vertex $w$ enters vertex $v$.
	\end{defn}

    \begin{ex}
    	Let's take a quick look at an tiny example.
    	
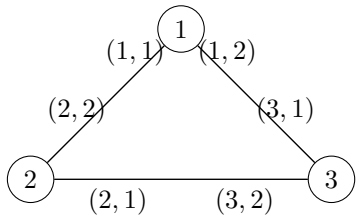
\[
    	\begin{tikzpicture}
    	\node[circle,draw](1) at (0,0){1};
    	\node[circle,draw](2) at (-2,-2){2};
    	\node[circle,draw](3) at (2,-2){3};
    	
    	\path[-,draw, thin] (1) edge node[near start,above]{$(1,1)$} (2); 
    	\path[-,draw, thin] (2) edge node[near start, above]{$(2,2)$} (1);
    	\path[-,draw, thin] (1) edge node[near start,above]{$(1,2)$} (3); 
    	\path[-,draw, thin] (3) edge node[near start, above]{$(3,1)$} (1);
    	\path[-,draw, thin] (2) edge node[near start,below]{$(2,1)$} (3); 
    	\path[-,draw, thin] (3) edge node[near start, below]{$(3,2)$} (2);
    	\end{tikzpicture}
    	\]
    	
    	The Rotation map can be given in a table as
    	\[
    	\begin{tabular}{l | l}
    	    $(v,i)$ & $(w,j)$\\
    	    \hline
    		(1,1) & (2,2)\\
    		(1,2) & (3,1)\\
    		(2,1) & (3,2)\\
    		(2,2) & (1,1)\\
    		(3,1) & (1,2)\\
    		(3,2) & (2,1)  
    	\end{tabular}
    	\]

    \end{ex}

    We can see a few things immediately.  First, the rotation map is an involution.  That is 
    \[
    Rot_G \circ Rot_G (v,i) = (v,i)
    \]

	Second, since this is an involution we can reduce this into a matrix.  In the above example we'll rewrite this as
	\[
	\begin{tabular}{l | c c |}
	$v\backslash i$  & 1 & 2\\
	  \hline
	1 & 2 & 3\\
	2 & 3 & 1\\
	3 & 1 & 2
	\end{tabular}
	\]

	That is we can reduce our rotation map by defining it as a matrix.
	\begin{defn}
		Given a regular graph $G$ with degree of regularity $d$ we redefine the rotation map as a matrix $Rot^{(G)} \in \Z^{|V|\times d}$
		by
		\begin{equation}
		Rot^{(G)}_{v,i} = w \label{rotation matrix}
		\end{equation}
		where $Rot_G(v,i) = (w,j)$ in the previous definition.  Since $Rot_G$ is an involution the information $j$ is superfluous.
	
	\end{defn}

	So far, we have only the definition of a rotation map, but nothing of consistency.  Consistent rotation maps are of paramount importance for this work so we consider
	\begin{defn}
		A rotation map $Rot_G$ is said to be \emph{consistent} if every vertex has $d$ distinct incoming labels.
	\end{defn}

    We can see that the rotation map in the above example is consistent.  In particular, since this is a 2-regular graph, we see that each vertex has two incoming edges with two distinct labels.  I.e.  The vertex $(1)$ has an incoming edge labeled1 coming from vertex $(3)$ and an incoming vertex labeled 2 coming from vertex $(2)$.

    In our new formulation, this is an easily checkable criterion.
    \begin{lem}
    	A rotation map given in matrix form is consistent if each column contains every element in $[d]$.
    \end{lem}

	\subsection{Reading a Rotation Map from the Adjacency Matrix}
	
	Given a graph $G$ and its corresponding adjacency matrix $A$, we can easily construct a rotation map $Rot^{(G)}$ in the style of equation \ref{rotation matrix}.
	
	This is given by the following algorithm.
	\begin{algorithm}
		\caption{Rotation map from the adjacency matrix; Julia/Octave style}
		\label{Adjacency to rotation}
		\begin{algorithmic}
			\State for $k$ in $1:|V|$
			\State $\hspace{0.5cm}$ counter = 1
			\State for $l$ in $1:|V|$
			\State $\hspace{0.5cm}$if $A[k,l] == 1$
			\State $\hspace{1cm} Rot^{(G)}[v,\text{counter}] = l$
			\State $\hspace{1cm}$ counter $+= 1$
			\State endif
			\State endfor
			\State endfor
		\end{algorithmic}	
	\end{algorithm}

	Notice that we have given this in syntax which is mostly compatible with julialang and Octave/MatLab since these are 1-indexed languages.
	
	We notice, however, that reading a rotation map in this way, will never produce a consistent rotation map. In particular the vertex labeled $(1)$ will show up $d$ times in the left most column and the vertex labeled $(|V|)$ will show up exactly $d$ times in the right most column.

	\subsection{Some Algebraic Properties of Cartesian Products}
	
	Given two graphs $G_1 = (V_1,E_1)$ and $G_2 = (V_2,E_2)$ which are both regular of degrees $d_1$ and $d_2$  and adjacency matrices $A_1$ and $A_2$ respectively their Cartesian product $G_1 \square G_2$ is defined by the following.
	
	\begin{defn}
		The Cartesian product $G_1\square G_2$ is the graph corresponding to the adjacency matrix
		\begin{equation}
		A_{1\square 2} = A_1 \otimes I_{|V_2|} + I_{|V_1|}\otimes A_2
		\end{equation}
	\end{defn}
	
	We get some immediate consequences using this definition.
	\begin{lem}
		Given two regular graphs $G_1 = (V_1,E_1)$ and $ G_2 = (V_2,E_2)$
		the Cartesian product $G_1\square G_2$ has the following properties
		\begin{itemize}
			\item[$\cdot$] $G_1 \square G_2$ has $|V_1|*|V_2|$ vertices.
			\item[$\cdot$] $G_1 \square G_2$ is a regular graph with degree of regularity $d_1 + d_2$.
			\item[$\cdot$] $G_1 \square G_2$ has $(|V_1|*|V_2|)(d_1+ d_2)/2$ edges.
			\item[$\cdot$] The spectrum of $G_1 \square G_2$ is given by
			\[
			\Lambda(G_1\square G_2) = \{\lambda + \mu   |  \lambda \in \Lambda(G_1), \mu\in \Lambda(G_2)\}
			\]
			where $\Lambda(G_i)$ are the eigenvalues of the adjacency matrix of $G_i$ taken with respective multiplicities.
		\end{itemize}	
	\end{lem}

	\subsection{The Main Method}
	
	The main algorithm is rather simple.  Following the lead of \cite{zig-zag} we use the notion of \emph{clouds} (or nuages).
	\begin{defn}
		Let $G = (V_1,E_1)$ and $H = (V_2,E_2)$ be two regular graphs with regularities $d_1,d_2$ respectively.  The product $G\square H$ is a graph with $V_1V_2$ vertices and degree $d_1 + d_2$.  We partition the  vertices of $G\square H$ into $V_2$ groups of $V_1$ vertices via the naive partition 
		\[
		\{v_1, ... , v_n\} = \{v_1,...v_{|G|}\},\{v_{|G|+1}, \dots v_{2|G|}\},\dots \{v_{1+ (|H|-1)|G|},\dots v_{|G||H|}\}
		\]
		each partition is called a \emph{cloud}.
	\end{defn}
	
	The idea is that each vertex of $H$ is replaced with a ``cloud" of the graph $G$.  Let's take an extremely simple example.  Let $G = C_3$ and $H = C_4$.  
	
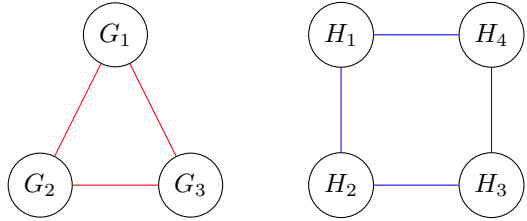
\[
	\begin{tikzpicture}
	\node[circle,draw](1) at (0,0){$G_1$};
	\node[circle,draw](2) at (-1,-2){$G_2$};
	\node[circle,draw](3) at (1,-2){$G_3$};

    \node[circle,draw](4) at (3,0){$H_1$};
    \node[circle,draw](5) at (3,-2){$H_2$};
    \node[circle,draw](6) at (5,-2){$H_3$};
    \node[circle,draw](7) at (5,0){$H_4$};
	
	\foreach \from/\to in {1/2,1/3,3/2} 
	\draw[red] (\from)--(\to);
	\foreach \from/\to in {4/5,5/6,6/7,7/4}
	\draw[blue] (\from)--(\to);
	\end{tikzpicture}
	\]
	
	The product will have four clouds of three vertices labeled below as $C_{i,j}$ that is,vertex $j$ in cloud $i$.
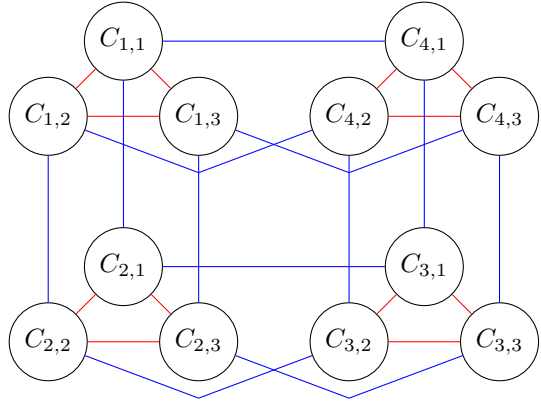
	\[
	\begin{tikzpicture}
	\node[circle, draw](1) at (0,0) {$C_{1,1}$};
	\node[circle, draw](2) at (-1,-1) {$C_{1,2}$};
	\node[circle, draw](3) at (1,-1) {$C_{1,3}$};
	\node[circle, draw](4) at (0,-3) {$C_{2,1}$};
	\node[circle, draw](5) at (-1,-4) {$C_{2,2}$};
	\node[circle, draw](6) at (1,-4) {$C_{2,3}$};
	\node[circle, draw](7) at (4,0) {$C_{4,1}$};
	\node[circle, draw](8) at (3,-1) {$C_{4,2}$};
	\node[circle, draw](9) at (5,-1) {$C_{4,3}$};
	\node[circle, draw](10) at (4,-3) {$C_{3,1}$};
	\node[circle, draw](11) at (3,-4) {$C_{3,2}$};
	\node[circle, draw](12) at (5,-4) {$C_{3,3}$};
	
	\foreach \from/\to in {1/2, 2/3, 3/1, 4/5,5/6,6/4,7/8,8/9,9/7,10/11,11/12,12/10} 
	\draw[red](\from)--(\to);
	
	\foreach \from/\to in {1/4,1/7, 7/10, 4/10, 2/5, 8/11, 3/6, 9/12}
	\draw[blue] (\from)--(\to);

	\draw[blue] (2)  --(1,-1.75) --  (8);
	\draw[blue] (3)  --(3,-1.75) -- (9);
	\draw[blue] (5)  --(1,-4.75) -- (11);
	\draw[blue] (6)  --(3,-4.75) -- (12);
	
	\end{tikzpicture}
	\]

	We will write our new rotation map in terms of the constituent rotation maps taking into account our idea of clouds.
	
	\begin{lem}
		Given two regular graphs $G$ and $H$ with consistent rotations maps $Rot^{(G)}$ and $Rot^{(H)}$ given in matrix form, with form the new rotation map $Rot^{(G\square H)}$ as a block matrix with $|V_H| \times 2$ blocks.  The blocks in the first column are of size $|V_G|\times d_G$ and blocks in the second column are of size $|V_G|\times d_H$. Denote these blocks as $B_{i,j}$.  We have two distinct formulas.
		\begin{enumerate}
			\item $B_{i,1} = Rot^{(G)} + (i-1)*|V_G|$. That is, each block in the first column is a copy of $Rot^{(G)}$ where we simply add a constant to each element.  
			\item $B_{i,2}$ row $j$, column $k = $ vertex $j$ of cloud $\ell$ where 
			\begin{equation}
			\ell = Rot^{(H)}_{i,k}
			\end{equation}
			This is slightly more complicated, but the basic idea is to connect the vertices in each cloud to corresponding vertices in the other clouds according to the rotation map given by $H$.
		\end{enumerate}
		 
	\end{lem}

	Algorithmically we have
	\begin{algorithm}
		\caption{Main Algorithm, Rotation Maps on Cartesian Products}
		\label{Rotation on Product}
		\begin{algorithmic}
			\State input: $Rot^{(G)}$, $Rot^{(H)}$ in matrix form.
			\State for $i$ in 1 to $|V_H|$
			\State $\hspace{0.5cm}$ cloud$[i] := \{(i-1)|V_G| + 1, \dots, i|V_G| \}$
			\State $\hspace{0.5cm}$ Block$[i,1] = Rot^{(G)} + (i-1)*ones(|V_G|,d_G)$
			\State endfor
			\State for $i$ in 1 to $|V_H|$
			\State for $j$ in 1 to $|V_G|$
			\State for $k$ in 1 to $d_H$
			\State $\hspace{0.5cm}\ell = Rot^{(H)}[i,k]$
			\State $\hspace{0.5cm}$ Block$[i,2][j,k] = j + (\ell-1)|V_G|$
			\State endfor
			\State endfor
			\State endfor
		\end{algorithmic}	
	\end{algorithm}

	\section{Families of Graphs with Easy Constructions}
	Here we give a few families of graphs with easily constructible consistent rotation maps.
	\begin{enumerate}
		\item Cyclic graphs on $n$-vertices:\\
		Label the vertices 1 to $n$ in a circle.  Since this is two-regular we have two columns.  Column 1 is the array $[2,3,\dots,n,1]$ and column 2 is the array $[n,1,2,\dots,n-1]$\\ 
		E.g. $G = C_5$ 
		\[ Rot^{(G)} = 
		\begin{bmatrix}
		2&5\\3&1\\4&2\\5&3\\1&4
		\end{bmatrix}
		\]

		\item Complete graphs on $n$ vertices:\\
		Again, we label the vertices 1 to $n$ in a circle.  The idea is to go around the edges clockwise or counterclockwise (or anticlockwise) labeling the edges 1 to $n-1$.
		This means the columns follow the same array patterns as the cyclic graph on $n$ vertices with the first row being $[2,3,\dots, n]$.
		E.g. $G = K_5$
		\[
		Rot^{(G)} = \begin{bmatrix}
		2&3&4&5\\
		3&4&5&1\\
		4&5&1&2\\
		5&1&2&3\\
		1&2&3&4
		\end{bmatrix}
		\]
		
		Notice that columns 1 and 5 are identical to columns 1 and 2 from the $C_5$.

		\item Complete bipartite graphs on $n,n$ vertices:\\
		Split the vertices into two columns, left and right.  The left column is labeled 1 to $n$.  The right column is labeled $n+1$ to $2n$.  For each vertex, the edge drawn horizontally is labeled 1, and then we proceed clockwise (or counter clockwise) as before.
		
		E.g. $G = K_{3,3}$
		\[
		Rot^{(G)} = \begin{bmatrix}
		4&5&6\\5&6&4\\6&4&5\\1&2&3\\2&3&1\\3&1&2
		\end{bmatrix}
		\]

		\item Generalized Petersen graphs:\\
		The general layout of a generalized Petersen graph is two concentric circles.
		of vertices labeled 1 to $n$ on the outer circle and $n+1$ to $2n$ on the inner circle.  The vertices $j$ and $n+j$ are connected by an edge.  On the outer circle, we connect vertices in a circle as above.  In the inner circle we may shift by any of $1$ to $n \text{div}2$ vertices. The graphs are label $GP(n, s)$ which is a graph on $2n$ vertices and the inner circle is shifted by $s$.  The Petersen graph is $GP(5,2)$.  The graphs $GP(n,1)$ are the cartesian products $C_n \times K_2$.
		Since these are all cubic graphs, we will set the edges $j$ with $n+j$ as `2' in both directions.  The first column will be identical to the cyclic graph on $n$ vertices and then a cycle around $n+1$ to $2n$ shifting each by $s$.  The third column is the remaining vertex. \\
		E.g. $G = GP(7,3)$ This is a 14 vertex graph where the inner circle makes a star in the order $[8,11,14,10,13,9,12]$ which gives us the rotation map
		\[
		Rot^{(G)} = \begin{bmatrix}
		2 & 8 & 7\\ 
		3 & 9 & 1\\
		4 & 10 & 2\\
		5 & 11 & 3\\
		6 & 12 & 4\\
		7 & 13 & 5\\
		1 & 14 & 6\\
		11 & 1 & 12\\
		12 & 2 & 13\\
		13 & 3 & 14\\
		14 & 4 & 8\\
		8 & 5 & 9\\
		9 & 6 & 10\\
		10 & 7 & 11\\
		\end{bmatrix}
		\]

	\end{enumerate}
	
	
	\section{Future Work}
	The current literature on rotation maps of regular graphs is thin.  It is an open question as to whether there is a constructive method for computing consistent rotation maps efficiently (ie not brute force).  The next problems which remain to be solved are if given a random regular graph, one can construct a consistent rotation map.  The utility from the perspective of the author is that once a consistent rotation map is known one can efficiently compute discrete time coined quantum walks on regular graphs.  This will allow those interested in DTQW to compute them on standard laptops. As of 2020 for even moderately sized graphs DTQW requires roughly $O(N^4)$ memory.  The ability of the user to compute DTQW will be in the ability of a modern computing language to store complex matrices.  The author makes no claim as to efficiently processing sparse matrices, however, even in the case of sparse matrices, a $d$-regular graph on $N$ vertices will produce a unitary walking matrix of size $Nd\times Nd$ versus $N^2\times N^2$ in other DTQW scenarios.  In particular a ``buckyball" which is 3-regular on 60 vertices can be easily computed in a $180\times 180$ complex matrix, which is extremely lightweight for a numerical engine such as MatLab/Octave, Julia, Numpy etc. Contrast this with a weighted coined walk in the style of Wong \cite{wong} which requires a $3600\times 3600$ matrix for implementation.    
	
	The next chapters in the rotation maps saga are to produce a heuristic solver for rotation maps on random regular graphs and on regular skeletons of manifolds.  With these constructions in hand, there are several hypothesized applications.
	\begin{itemize}
		\item[(a)] In Materials science, there are several important molecules which can be constructed as regular graphs, namely carbon nanostructures which are 3-regular (cubic) graphs.  The first hypothesis is that rigorous calculations of physical properties can be made by properly constructing DTQW.  This follows the lead of Bellisard \cite{Bellisard}et el, but is a different technique entirely.
		\item[(b)] Calculation of topological properties of manifolds.  In particular the hypothesis is that with a specific initial state and initial coin, one can, with high probability compute $H_1(M)$.
		\item[(c)] ``More random" random sampling.  In particular optimization by random sampling is a deeply studied subject (cf Metropolis) and has new techniques (MCMCMC) \cite{mc3}.  The hypothesis is that one can produce a sequence of random chains in an MCMCMC and land on a  desired distribution ``faster."
	\end{itemize}

	\section*{An Example: The Skeleton of a Torus}
	In this section we'll construct two examples from the ground up.  First,let's look at the torus.  We're constructing the skeleton of the torus, which can be thought of $S^1\times S^1$ which we discretize as $C_n \times C_m$.  Note that $n$ need not equal $m$.  There are some interesting dynamics in discrete time quantum walks on tori which depend on the ratio $\frac{n}{m}$.  We will give the general frame work, build up a specific example and allow the interested reader to construct multiple tori at one's leisure.  
	
	\begin{ex}
		Let's consider the torus of $C_6 \times C_4$.  This will be a 4 regular graph with 24 vertices. Our two rotation maps as from the construction above (\S 3)
		\begin{equation}
		Rot^{(C_6)} = \begin{bmatrix}
		2&6\\3&1\\4&2\\5&3\\6&4\\1&5
		\end{bmatrix}
		\text{ and }
		Rot^{(C_4)} = \begin{bmatrix}
		2&4\\3&1\\4&2\\1&3
		\end{bmatrix}
		\end{equation}
		
	We have 4 clouds 
	\begin{enumerate}
		\item[Cloud 1:] [1, 2,3, 4, 5, 6]
		\item[Cloud 2:] [7, 8, 9, 10, 11, 12]
		\item[Cloud 3:] [13, 14, 15, 16, 17, 18]
		\item[Cloud 4:] [19, 20, 21, 22, 23, 24]
	\end{enumerate}
	and 8 blocks to construct.  As mentioned in algorithm 2, the block $B_{i,1}$ are simple to construct.
	
	\begin{equation}
	B_{1,1} = \begin{bmatrix}
	2&6\\3&1\\4&2\\5&3\\6&4\\1&5
	\end{bmatrix} \text{ } 
	B_{2,1} = \begin{bmatrix}
	8&12\\9&7\\10&8\\11&9\\12&10\\7&11
	\end{bmatrix} \text{ }
	B_{3,1} = \begin{bmatrix}
	14&18\\15&13\\16&14\\17&15\\18&16\\13&17
	\end{bmatrix} \text{ }
	B_{4,1} = \begin{bmatrix}
	20&24\\21&19\\22&20\\23&21\\24&22\\19&23
	\end{bmatrix} 
	\end{equation}
	
	The second column of blocks are slightly trickier.  Let's look closely at $B_{1,2}$.
	By our lemma we will be concerned with row 1 of $Rot^{(C_4)}$ ie $[2,4]$.
	Then $B_{1,2}[j,k] = j+(\ell-1)|V_G|$  where $\ell = $ element $k$ of $[2,4]$ and $|V_G|=6$.
	
	Therefore we have
	\begin{eqnarray}
	B_{1,2}[j,k] & = & j + 6(\ell-1)\\
	\nonumber B_{1,2}[1,1] & = & 1 + 6(2-1)\\
	\nonumber B_{1,2}[1,2] & = & 1 + 6(4-1)\\
	\nonumber B_{1,2}[2,1] & = & 2 + 6(2-1)\\
	\nonumber B_{1,2}[2,2] & = & 2 + 6(4-1)\\
	\nonumber B_{1,2}[3,1] & = & 3 + 6(2-1)\\
	\nonumber B_{1,2}[3,2] & = & 3 + 6(4-1)\\
	\nonumber B_{1,2}[4,1] & = & 4 + 6(2-1)\\
	\nonumber B_{1,2}[4,2] & = & 4 + 6(4-1)\\
	\nonumber B_{1,2}[5,1] & = & 5 + 6(2-1)\\
	\nonumber B_{1,2}[5,2] & = & 5 + 6(4-1)\\
	\nonumber B_{1,2}[6,1] & = & 6 + 6(2-1)\\
	\nonumber B_{1,2}[6,2] & = & 6 + 6(4-1)
	\end{eqnarray}
	
	Which reveals
	\begin{equation}
	B_{1,2} = \begin{bmatrix}
	7&19\\
	8&20\\
	9&21\\
	10&22\\
	11&23\\
	12&24
	\end{bmatrix}
	\end{equation}
	
	Repeating, we get the final three blocks
	\begin{equation}
	B_{2,2} = \begin{bmatrix}
	13&1\\14&2\\15&3\\16&4\\17&5\\18&6
	\end{bmatrix} \text{ } 
	B_{3,2} = \begin{bmatrix}
	19&7\\20&8\\21&9\\22&10\\23&11\\24&12
	\end{bmatrix} \text{ } 
	B_{4,2} = \begin{bmatrix}
	1&13\\2&14\\3&15\\4&16\\5&16\\6&18
	\end{bmatrix}  
	\end{equation}
		
	\end{ex}


\begin{thebibliography}{80}
        \bibitem[A]{A} Ambainis, A. \emph{Quantum Walks and their Algorithmic Applications}, https://arxiv.org/abs/quant-ph/0403120
        
        \bibitem[AKR]{AKR} Ambainis, A, Kempe, J., Rivosh, A., \emph{Coins Make Quantum Walks Faster}, https://arxiv.org/abs/quant-ph/0402107	
	
	    \bibitem[B]{Bellisard} Bellisard, J., \emph{The Noncommutative Geometry of Aperiodic Solids}, "Geometric and Topological Methods for Quantum Field Theory", (Villa de Leyva, 2001), pp. 86-156,
	    World Sci. Publishing, River Edge, NJ, (2003).
	    
	    \bibitem[G]{mc3} Geyer, C. J. \emph{Markov chain Monte Carlo maximum likelihood}, Computing Science and Statistics: Proc. 23rd Symp. Interface, 156–163, 1991
	    
	    \bibitem[K]{K} Kendon, V., \emph{Quantum Walks on General Graphs}, https://arxiv.org/abs/quant-ph/0306140
	
		\bibitem[RVW]{zig-zag}Reingold, O.; Vadhan, S.; Widgerson, A. (2000), \emph{Entropy waves, the zig-zag graph product, and new constant-degree expanders and extractors}, 41st Annual Symposium on Foundations of Computer Science: 3–13, arXiv:math/0406038, doi:10.1109/SFCS.2000.892006, ISBN 978-0-7695-0850-4
		
		\bibitem[W]{wong} Wong, T. \emph{Coined Quantum Walks on Weighted Graphs}, 2017 J. Phys. A: Math. Theor. 50 475301
\end{thebibliography}
\end{document}